\documentclass[11pt]{amsart}
\usepackage{amsmath, amssymb, amsthm,latexsym}
\usepackage{graphicx,color}

\newcommand\N{{\mathbb N}}

\newcommand\R{{\mathbb R}}

\newcommand\dee{\partial}

\newcommand\id{\operatorname{id}}

\renewcommand\:{\colon}

\newcommand\no{\noindent}

\newtheorem{theorem}{Theorem}[section]

\newtheorem{proposition}[theorem]{Proposition}
\newtheorem{corollary}[theorem]{Corollary}

\newtheorem{lemma}[theorem]{Lemma}

\theoremstyle{definition}

\title[Hausdorff dimension of three-period orbits]{Hausdorff dimension of three-period orbits in Birkhoff billiards}

\author{Sergei Merenkov}
\address{Sergei Merenkov\\Department of Mathematics\\
University of Illinois\\ 1409 W~Green Street\\ Urbana, IL 61801
\\USA} \email{merenkov@illinois.edu}
\thanks{S.M.\ was supported by NSF grant DMS-1001144.}

\author{Vadim Zharnitsky}
\address{Vadim Zharnitsky\\Department of Mathematics\\
University of Illinois\\ 1409 W~Green Street\\ Urbana, IL 61801
\\USA} \email{vz@math.uiuc.edu}

\begin{document}

\begin{abstract}
We prove that the Hausdorff dimension of the set of  three-period orbits 
in classical billiards is at most one. 
Moreover, if the set of  three-period orbits has Hausdorff dimension one, then it has a tangent line at almost every point. 

\end{abstract}

\maketitle

\section{Introduction}
\no
The goal of this note is to provide a sharp estimate on the dimension of the 
set $P^3$ of three-period orbits in classical billiards. The question about Lebesgue
 measure of the set of periodic  orbits was raised by Ivrii~\cite{ivrii}, in connection
 with spectral geometry problems. 
The conjecture that periodic orbits have zero measure is still open, but 
several  special cases have been resolved. The first result was due to Rychlik \cite{rychlik} who proved that the set of three-period orbits has zero Lebesgue measure. 
More recently a proof for the four-period case has been announced in \cite{glku}.

Consider a not necessarily convex billiard boundary and define the billiard map $T$
in the standard way. The boundary of length $l$ is parametrized with a natural parameter $t$. 
For an outgoing trajectory, let $\theta$ be the angle with the tangent. Then the billiard map $T(t_1,\theta_1) = (t_2,\theta_2)$ is defined 
for any segment that is transversal to the boundary at both ends. We disregard those that are not, since 
they have zero measure. The billiard map is defined on an open annulus $U= S^1\times I$, where $S^1$ is 
the boundary of length $l$  and $I=(0,\pi).$  The billiard map preserves the measure $\mu = \sin(\theta) d\theta \wedge dt$, see, e.g.,  \cite{tabachnikov}. 

Rychlik's proof was simplified  in \cite{stojanov}, \cite{vorobets}, and \cite{wojtk}, see also a survey article \cite{gutkin_s}. 
In \cite{wojtk},  Wojtkovski proved Rychlik's theorem  using Jacobi fields. He first proved that 
 if there is a neighborhood containing only three-period orbits, then the identity
\[
k(t) \cdot L(t,\theta) = 2 \sin^3(\theta)
\]
must hold, 
where $k(t)$ is the boundary curvature and $L(t,\theta)$  is the length of the three-period orbit $(t,\theta)$. 
By Fermat's principle, the length is extremized by the actual orbits, and thus  $L(t,\theta)$ is constant
for any continuous family of orbits. 
The contradiction is easily derived by observing that, by varying $\theta$, 
the right-hand side will change while the left-hand side will be constant. Wojtkovski, following Rychlik 
\cite{rychlik}, 
extended the contradiction to the case when the set $P^3$ has positive Lebesgue measure. The argument relies 
on the presence of a Lebesgue density point in whose neighborhood one can differentiate 
the above identity in $\theta$, arriving at a contradiction again. 

A natural question arises as to how optimal this result is. Clearly, periodic orbits can form 
a one-parameter  family, e.g., in a billiard with circular boundary there is a one dimensional set of three-period orbits, so the  Hausdorff dimension can be equal to one. But can this set have the Hausdorff dimension between one and two?
 Note that for the two-period orbits the answer is clearly no, since such orbits must be normal 
to the boundary. Thus, $\theta = \pi/2$ and the set is  confined to the one-dimensional  line segment.

Our main result states that the similar sharp bound holds for the three-period orbits.

\begin{theorem}\label{T:main}
If the billiard boundary is $\mathcal C^3$, then the Hausdorff dimension of 
the set of three-period orbits $P^3$ is not more than one, i.e.,
\[
{\mathcal H}^s (P^3) = 0 \,\,\, {\rm if} \,\,\, s > 1,
\]
where $\mathcal H^s$ denotes the Hausdorff $s$-dimensional measure.
\end{theorem}

In the case when the Hausdorff dimension is equal to one, 
there is more information about the structure  of the set $P^3$, exhibited by the following theorem.
See Section~\ref{S:MT} for the definitions. 

\begin{theorem}\label{T:1reg}
Assume the billiard boundary is $\mathcal C^3$ and 
${\rm dim}_{\mathcal H^1}(P^3)=1$.
Then the set $P^3$ has a tangent line at  $\mathcal H^1$-almost every point $p\in P^3$.

\end{theorem}

Note that in the case of two-period orbits the analogous statement holds for the above reason: $P^2$ is  confined to the line segment $\theta = \pi/2$.

{
The paper is organized as follows. 
Section~\ref{S:Diff} contains an elementary lemma on existence of a sequence converging to an accumulation point of a given set and that is asymptotic to a certain direction. It also provides an elementary formula for a directional derivative of a map via a sequence of points accumulating asymptotically in the given direction. 
In Section~\ref{S:MT} we define and develop the relevant notions from measure theory, mostly based on \cite{kF86}. In the last section we apply those results to prove main theorems on the structure of the set of three-period orbits.}

\section{Derivatives along sets}\label{S:Diff}
\no
For two points $p$ and $q$ in $\R^2$, we denote by $|p-q|$ the Euclidean distance between them. 
If $E$ and $F$ are two subsets of $\R^2$, we denote by ${\rm dist}(E,F)$ the distance between them. We also denote by ${\rm diam}(E)$ the Euclidean diameter of a set $E\subseteq\R^2$.
Finally, let $B(p,r)$ denote the closed disc in $\R^2$ centered at $p$ of radius $r\geq0$.

For $\gamma\in\R$ and $\rho, \eta>0$, we consider
the {circular sector}
$$
W_p({\gamma; \rho, \eta})=\{q=p+te^{i\theta}\colon 0\leq t\leq\rho,|\theta-\gamma|\leq\eta\}.
$$

\begin{lemma}\label{L:TL}
Let $E$ be a subset of the disc $B(p,r)\subseteq\R^2$. Assume that $F$ is a closed subset of the boundary circle $\dee B(p,r)$ such that $p$ is an accumulation point for $E\cap C_F$, where $C_F$ is the cone over $F$, i.e.,
$$
C_F=\{q\in B(p,r)\colon q=p+te^{i\theta},\ re^{i\theta}\in F,\ t\in[0,r]\}.
$$ 
Then there exist a ray $R$, emanating from $p$ and passing through a point in $F$, and a sequence $(p_k)$ with $p_k\in E\setminus\{p\}$, $\lim p_k=p$, such that $(p_k)$ is asymptotic to $R$, i.e.,
$$
\lim_{k\to\infty}\frac{{\rm dist}(p_k,R)}{|p_k-p|}=0.
$$
\end{lemma}
\no
\emph{Proof.}
For every $n\in\N$, we consider the open cover 
$$
\mathcal G_n=\{U_n(\gamma)\colon re^{i\gamma}\in F\}
$$ 
of $F$, where 
$$
{U}_n(\gamma)=\{re^{i\theta}\colon \theta\in(\gamma-1/2^n,\gamma+1/2^n)\}.
$$ 
Since $F$ is compact, $\mathcal G_n$ contains a finite subcover, denoted $\mathcal G_n(F)$.   Let $\mathcal G_n(F)=\{U_n(\gamma_n^1),\dots, U_n(\gamma_n^{N(n)})\}$.

We assumed that $p$ is an accumulation point for $E\cap C_F$. Therefore, for every $n\in\N$, there exists $\gamma_n\in\{\gamma_n^1,\dots, \gamma_n^{N(n)}\}$, such that $p$ is an accumulation point for 
$$
E\cap W_p(\gamma_n;1/n,1/2^n). 
$$
Moreover, we may assume that for each $n\in\N$, we have $re^{i\gamma_{n+1}}\in \overline{U_n(\gamma_n)}$, the closure of $U_n(\gamma_n)$.
Indeed, $F_n=F\cap\overline{U_n(\gamma_n)}$ is a compact set such that $p$ is an accumulation point for $E\cap C_{F_n}$. Thus the sequence $(re^{i\gamma_n})$ converges to $re^{i\gamma_\infty}$. Note that  $re^{i\gamma_\infty}\in F$.

Let $R$ be the ray emanating from $p$ and passing through the point $re^{i\gamma_\infty}$. From the above, we know that there exists 
a sequence $(p_n),\ p_n\in E\setminus\{p\}$, such that 
$$
p_n\in W_p(\gamma_{n};1/n, 1/2^{n}).
$$
Since $\lim 1/n=0$, we conclude that $\lim p_n=p$. Finally,
$$
\frac{{\rm dist}(p_n,R)}{|p_n-p|}\leq\sin\left(\frac{1}{2^{n}}+\sum_{j=0}^\infty\frac1{2^{n+j}}\right)=\sin\left(\frac3{2^{n}}\right)\to0,\quad k\to\infty. 
\qed
$$

The following lemma is elementary.
\begin{lemma}\label{L:DD}
Let $p$ be an arbitrary point in $\R^2$ and let $U\subseteq\R^2$ be a neighborhood of $p$.
Suppose that $F\colon U\to\R$ or $\R^2$ is a $\mathcal C^1$-differentiable map
in $U$. Let $\bf{v}$ be a unit vector in $\R^2$ and $R$ be a ray emanating from $p$ in the direction of $\bf{v}$. We assume that $(p_k)$ is a sequence in $U\setminus\{p\}$, such that $\lim_{k\to\infty}p_k=p$ and $(p_k)$ is asymptotic to $R$ in the sense of Lemma~\ref{L:TL}. Then for the directional derivative $F_{\bf{v}}(p)=D_pF({\bf v})$, where $D_pF$ denotes the differential of $F$ at $p$, we have
$$
F_{\bf{v}}(p)=\lim_{k\to\infty}\frac{F(p_k)-F(p)}{|p_k-p|}.
$$ 
\end{lemma}
\no
\emph{Proof.}
For $k\in \N$, let $q_k$ be the projection of $p_k$ onto $R$. Then
$$
F_{\bf{v}}(p)=\lim_{k\to\infty}\frac{F(q_k)-F(p)}{|q_k-p|}.
$$
On the other hand,
\begin{equation}\label{E:DD}
\frac{F(p_k)-F(p)}{|p_k-p|}=\frac{F(q_k)-F(p)}{|q_k-p|(1+o(1))}+\frac{F(p_k)-F(q_k)}{|p_k-p|},\quad k\to\infty,
\end{equation}
because $(p_k)$ is asymptotic to $R$. Also, since $F\in\mathcal C^1$ in $U$, it is Lipschitz in a neighborhood of $p$. Therefore there exists $L>0$ such that 
$$
|F(p_k)-F(q_k)|\leq L|p_k-q_k|=L\, {\rm dist}(p_k,R),
$$
for all $k$ large enough. By passing to the limit in~(\ref{E:DD}) and using the assumption that $(p_k)$ is asymptotic to $R$, we obtain the desired formula for the directional derivative $F_{\bf{v}}(p)$.
\qed

\section{Measure theory}\label{S:MT}
\no
Let $E\subseteq\R^2$ and $0\leq s<\infty$. For $\delta>0$ we define
$$
\mathcal H^s_\delta(E)=\inf\sum_{i\in I} {\rm diam}^s(U_i),
$$
where the infimum is taken over all covers $\{U_i\: i\in I\}$ of $E$ with $0<{\rm diam}(U_i)\leq\delta,\ i\in I$.  The \emph{Hausdorff} $s$-\emph{dimensional outer measure} of $E$ is
$$
\mathcal H^s(E)=\lim_{\delta\to0+}\mathcal H^s_\delta(E).
$$
The Hausdorff $s$-dimensional outer measure is indeed an \emph{outer measure}, i.e., $\mathcal H^s(\emptyset)=0$, $\mathcal H^s$ is monotone and subadditive. It is also a \emph{metric outer measure}, i.e., if $E$ and $F$ are two subsets of $\R^2$ with ${\rm dist}(E,F)>0$, then
$$
\mathcal H^s(E\cup F)=\mathcal H^s(E)+\mathcal H^s(F).
$$
Therefore we have the following consequence.
\begin{lemma}\label{L:MOM}\cite[Theorem~1.5]{kF86}
All Borel sets are $\mathcal H^s$-measurable.
\end{lemma}

It is easy to see that for any $E\subseteq\R^2$, the function $\mathcal H^s(E)$ is non-increasing in $s$. Furthermore, there exists a unique value $s_0$ such that $\mathcal H^s(E)=\infty$ for $0\leq s<s_0$ and $\mathcal H^s(E)=0$ is $s> s_0$. The value $s_0$ is called the \emph{Hausdorff dimension} of $E$, and it is denoted ${\rm dim}_\mathcal{H}(E)$.
An $s$-\emph{set} $E\subseteq\R^2$ is an $\mathcal H^s$-measurable set with $0<\mathcal H^s(E)<\infty$. It is immediate that if $E$ is an $s$-set, its Hausdorff dimension is $s$.

\begin{theorem}\label{T:SS}\cite[Theorem~5.4(a)]{kF86}
Let $E\subseteq\R^2$ be a closed set with $\mathcal H^s(E)=\infty$. Then for every $c>0$, there exists a compact subset $F\subseteq E$ such that $\mathcal H^s(F)=c$. In particular, $F$ is an $s$-set.
\end{theorem}

If $E$ is a $\mathcal H^s$-measurable set in $\R^2$, we say that $F\subseteq E$ is a \emph{full measure subset} if $\mathcal H^s(E\setminus F)=0$.
Also, if $E$ is a $\mathcal H^s$-measurable set, we say that a property holds at $\mathcal H^s$-\emph{almost all} points of $E$, or at $\mathcal H^s$-\emph{almost every} point of $E$, if there exists a full measure subset $F\subseteq E$ such that the property holds at every point of $F$.

If $E\subseteq\R^2$ is a $\mathcal H^s$-measurable set and $p\in\R^2$, the \emph{lower} and \emph{upper density} of $E$ at $p$ are defined as
$$
D_s(E,p)=\liminf_{r\to0+}\frac{\mathcal H^s(E\cap B(p,r))}{(2r)^s},\ 
D^{s}(E,p)=\limsup_{r\to0+}\frac{\mathcal H^s(E\cap B(p,r))}{(2r)^s},
$$ 
respectively.
If $E$ is an $s$-set, a  point $p\in E$ such that $D_s(E,p)=D^s(E,p)=1$ is called a \emph{regular} point of $E$. Otherwise $p$ is called an \emph{irregular} point. 

\begin{lemma}\label{L:UD}\cite[Corollary~2.5]{kF86}
If $E$ is an $s$-set in $\R^2$, then
$$
\frac1{2^s}\leq D^s(E,p)\leq 1
$$
at $\mathcal H^s$-almost all $p$ in $E$.
\end{lemma}



For $\gamma\in\R$ and $\eta>0$, we also define the \emph{upper angular density} of $E$ at $p$ as  
$$
D^s(E,p,\gamma,\eta)=\limsup_{r\to0+}\frac{\mathcal H^s(E\cap{W_p(\gamma;r,\eta)})}{(2r)^s}.
$$

A $\mathcal H^s$-measurable set $E\subseteq\R^2$ has a \emph{tangent line} at $p$ in direction $\gamma\in \R$ if $D^s(E,p)>0$ and if for every $\eta>0$,
$$
\lim_{r\to0+}\frac{\mathcal H^s ( E\cap(B(p,r)\setminus({W_p(\gamma;r,\eta)}\cup{W_p(\gamma+\pi;r,\eta)})))}{r^s}=0.
$$
We recall that $W_p({\gamma; \rho, \eta})$ is the circular sector
$$\{q=p+te^{i\theta}\colon 0\leq t\leq\rho,\ \gamma-\eta\leq\theta\leq \gamma+\eta\}.
$$

\begin{lemma}\label{L:AUDB}\cite[Lemma~4.5]{kF86}
If $1<s<2$ and $E$ is an $s$-set in $\R^2$, then for $\mathcal H^s$-almost all $p\in E$ we have 
$$
D^s(E,p,\gamma,\eta)\leq 4\cdot 10^s\eta^{s-1}
$$ 
for all $\gamma\in\R$ and $\eta\leq\pi/2$.
\end{lemma}

\begin{corollary}\label{C:TL}
If $1<s<2$ and $E$ is an $s$-set in $\R^2$, then at  $\mathcal H^s$-almost all points of $E$ no tangent line exists. 
\end{corollary}
\no
\emph{Proof.}
The sub-additivity of $\mathcal H^s$ gives 
\begin{equation}\label{E:SA}
\begin{aligned}
\mathcal H^s(E\cap B(p,r))&\leq \mathcal H^s(E\cap{W_p(\gamma;r,\eta)})+\mathcal H^s(E\cap{W_p(\gamma+\pi;r,\eta)})\\
&+\mathcal H^s(E\cap(B(p,r)\setminus({W_p(\gamma;r,\eta)}\cup{W_p(\gamma+\pi;r,\eta)})))
\end{aligned}
\end{equation}
for all $\gamma\in \R$ and $\eta>0$.

Suppose that $E$ has a tangent line $\gamma$ at $p$. Then, by dividing both sides of~(\ref{E:SA}) by $(2r)^s$ and taking $\limsup$ as $r\to0+$, we conclude that
$$
D^s(E,p)\leq D^s(E,p,\gamma,\eta)+D^s(E,p,\gamma+\pi,\eta)
$$
for all $\gamma\in \R$ and $\eta>0$.
By applying Lemma~\ref{L:AUDB}, we further obtain
$$
D^s(E,p)\leq 8\cdot 10^s\eta^{s-1},
$$
for all $\eta\leq\pi/2$ and $\mathcal H^s$-almost all $p\in E$ such that $E$ has a tangent line at $p$. Since $s>1$, this gives $D^s(E,p)=0$. An application of Lemma~\ref{L:UD} concludes the proof.
\qed


\begin{corollary}\label{C:UTL}
If $1<s<2$ and $E$ is an $s$-set in $\R^2$, then at $\mathcal H^s$-almost all points $p\in E$ we have
$$
0<\limsup_{r\to0+}\frac{\mathcal H^s(E\cap(B(p,r)\setminus({W_p(\gamma;r,\eta)}\cup{W_p(\gamma+\pi;r,\eta)})))}{r^s}\leq 2^s
$$
for all $\gamma\in\R$ and some $\eta>0$.
\end{corollary}
\no
\emph{Proof.}
The right-hand side inequality follows from the monotonicity of $\mathcal H^s$ and Lemma~\ref{L:UD}. The left-hand side inequality holds for every $p\in E$ such that $D^s(E,p)>0$ and $E$ has no tangent line at $p$. Lemma~\ref{L:UD} and Corollary~\ref{C:TL} imply that the set of such points has full $\mathcal H^s$-measure.  
\qed

\begin{theorem}\label{T:TATL}
Let $E\subseteq\R^2$ be an $s$-set, $1<s<2$. Then for $\mathcal H^s$-almost every point $p$ of $E$ there exist  two rays $R_1$ and $R_2$, emanating from $p$ and not contained in the same line, and two sequences $(p_k)$ and $(q_k)$ in $E\setminus\{p\}$, such that $\lim p_k=\lim q_k=p$, the sequence $(p_k)$ is asymptotic to $R_1$, and the sequence $(q_k)$ is asymptotic to $R_2$.
\end{theorem}
\no
\emph{Proof.}
Let $p\in E$ be an arbitrary point satisfying the two inequalities of Corollary~\ref{C:UTL} for all $\gamma\in\R$ and some $\eta>0$. That corollary states that the set of such points has full $\mathcal H^s$-measure.

Let $\gamma\in\R$ be arbitrary. Then there exists $\eta>0$ such that 
$$
\limsup_{r\to0+}\frac{\mathcal H^s(E\cap(B(p,r)\setminus({W_p(\gamma;r,\eta)}\cup{W_p(\gamma+\pi;r,\eta)})))}{r^s}>0.
$$
We may assume that $\eta<\pi/2$.
In particular, for a fixed $r>0$, $p$ is an accumulation point for 
$$
E\cap (W_p(\gamma+\pi/2; r,\pi/2-\eta)\cup W_p(\gamma-\pi/2; r,\pi/2-\eta)).
$$
By Lemma~\ref{L:TL}, there exist a ray $R_1$, emanating from $p$ and passing through a point in $$
F=\dee B(p,r)\cap (W_p(\gamma+\pi/2; r,\pi/2-\eta)\cup W_p(\gamma-\pi/2; r,\pi/2-\eta)), 
$$
and a sequence $(p_k)$ with $p_k\in E\setminus\{p\},\ \lim p_k=p$, such that $(p_k)$ is asymptotic to $R_1$. Clearly $R_1$ is not contained in the line through $p$ and $re^{i\gamma}$.

Let $\gamma'\in \R$ be such that $R_1$ contains the point $re^{i\gamma'}$. Then we can apply the same argument as above for $\gamma'$ in place of $\gamma$ to produce a ray $R_2$ and a sequence $(q_k)$ 
as in the statement of the theorem. The theorem follows.  
\qed

\begin{theorem}\label{T:SMD}
Let $U\subseteq\R^2$ be an open set and $F\: U\to\R^2$ be a $\mathcal C^1$-differentiable map. Let $E\subseteq U$ be a closed subset that consists of all fixed points of $F$. We assume that the Hausdorff dimension of $E$ is $s,\ 1<s\leq 2$. Then for every $s',\ 1<s'<s$, there exists an $s'$-set $E'\subseteq E$, such that $D_pF=\id$ for every $p\in E'$. 
\end{theorem}
\no
\emph{Proof.}
Let $s', 1<s'<s$ be arbitrary. Then $\mathcal H^{s'}(E)=\infty$, and therefore Theorem~\ref{T:SS} implies that there exists an $s'$-set $E_1\subseteq E$. Since $1<s'<2$, Theorem~\ref{T:TATL} yields the existence of a full measure set $E'\subseteq E_1$ with the following property.  For each point  $p\in E'$ there are two rays $R_1$ and $R_2$, emanating from $p$ and not contained in the same line, and two sequences $(p_k)$ and $(q_k)$ in $E'\setminus\{p\}$ with $\lim p_k=\lim q_k=p$, $(p_k)$ is asymptotic to $R_1$, and $(q_k)$ is  asymptotic to $R_2$.

Let $\bf{v_1}$ and $\bf{v_2}$ be the unit vectors that give the directions of $R_1$ and $R_2$, respectively. Note that $\bf{v_1}$ and $\bf{v_2}$ are linearly independent.
According to Lemma~\ref{L:DD}, for every $p\in E'$,
$$
F_{\bf{v_1}}(p)=\lim_{k\to\infty}\frac{F(p_k)-F(p)}{|p_k-p|}=\lim_{k\to\infty}\frac{p_k-p}{|p_k-p|}=\bf{v_1},
$$
since $E'\subseteq E$, and $(p_k)$ is asymptotic to $R_1$. Likewise, $F_{\bf{v_2}}(p)=\bf{v_2}$.

Since $\bf{v_1}$ and $\bf{v_2}$ form a basis in $\R^2$, for every $\bf{v}$ in $\R^2$ we have ${\bf{v}}=c_1{\bf{v_1}}+c_2{\bf{v_2}}$, where $c_1, c_2\in\R$. Thus, for every $p\in E'$, 
$$
D_pF({\bf{v}})=c_1D_pF({\bf{v_1}})+c_2D_pF({\bf{v_2}})= c_1F_{{\bf{v_1}}}(p)+c_2F_{{\bf{v_2}}}(p)=c_1{\bf{v_1}}+c_2{\bf{v_2}}=\bf{v},
$$
i.e., $D_pF=\id$.
\qed

Using a similar argument, one can prove the following result.
\begin{theorem}\label{T:SDF}
Let $U\subseteq\R^2$ be an open set and $f,g\: U\to\R$ be two $\mathcal C^1$-differentiable functions. We assume that $E\subseteq U$ is an $s$-set for $1<s<2$ and $f=g$ on $E$. Then there exists a full $\mathcal H^s$-measure subset $E'\subseteq E$ with the following property. For each $p\in E'$ there are two linearly independent unit vectors $\bf{v_1}$ and $\bf{v_2}$, such that
$f_{\bf{v_j}}(p)=g_{\bf{v_j}}(p),\ j=1,2$.  
\end{theorem}
\no
\emph{Proof.}
Let $E'\subseteq E$ be the full measure set that comes from  Theorem~\ref{T:TATL}, and let $p\in E'$ be arbitrary. Let $R_1$ and $R_2$ be two rays, emanating from $p$ and that are not contained in the same line. Let ${\bf{v_1}}$ and ${\bf{v_2}}$ be the linearly independent unit vectors that give the directions of $R_1$ and $R_2$, respectively. 
Let $R_j,\ j=1$ or 2, be one of the rays, and $(p_k),\ p_k\in E'\setminus\{p\},\ \lim p_k=p$, be an asymptotic sequence to $R_j$. 
Then, by applying Lemma~\ref{L:DD}, we get
$$
f_{\bf{v_j}}(p)=\lim_{k\to\infty}\frac{f(p_k)-f(p)}{|p_k-p|}=\frac{g(p_k)-g(p)}{|p_k-p|}=g_{\bf{v_j}}(p),
$$
as desired.
\qed

We note that the main ingredient of the proofs of Theorems~\ref{T:SMD} and~\ref{T:SDF} is the existence of sets of points at which no tangent line exists. Therefore we immediately have the following, more general, results.






\begin{theorem}\label{T:Diff}
Let $U\subseteq\R^2$ be an open set and $F\:U\to\R^2$ be a $\mathcal C^1$-differentiable map. Further, let $E\subseteq U$ be a closed subset that consists of all fixed points of $F$. 
Then $D_pF=\id$ at every $p\in E$ such that $E$ has no tangent line at $p$.
\end{theorem}

\begin{theorem}\label{T:Dir}
Let $U\subseteq\R^2$ be an open set and $f,g\: U\to \R$ be $\mathcal C^1$-differentiable functions. Then for every $p\in E$ such that $E$ does not have a tangent line at $p$, there are two linearly independent unit vectors $\bf{v_1}$ and $\bf{v_2}$, such that
$f_{\bf{v_j}}(p)=g_{\bf{v_j}}(p),\ j=1,2$.  
\end{theorem}

\section{Proofs of the main results}\label{S:Proofs}
\no
Assume that $ {\mathcal H}^s (P^3)\neq 0 $ for some $s\in (1,2]$. 
On the set $P^3 \subseteq U$ we have $T^3(p)=p$.
We may assume that the annulus $U=S^1\times I$ is smoothly embedded in the plane. We still denote the local coordinates by $(t,\theta)$.
Since the billiard boundary is $\mathcal C^3$, an elementary argument gives that the map $T^3$ is $\mathcal C^2$-differentiable, in particular $D_pT^3$ exists at each $p\in U$. We need the following result. 
\begin{proposition}\label{P:2func}\cite[Section 4]{wojtk}
If $T^3(t,\theta)=(t,\theta)$ and  $D_{(t,\theta)}T^3=\id$, then 
\begin{equation}\label{E:Main}
k(t) L(t,\theta) = 2 \sin^3(\theta),
\end{equation}
where $k(t)$ is the  boundary curvature at $t$ and $L(t,\theta)$ is the length of the three-period orbit $(t,\theta)$. 
\end{proposition}

\begin{figure}
[htbp]
\begin{center}
\includegraphics[height=50mm]{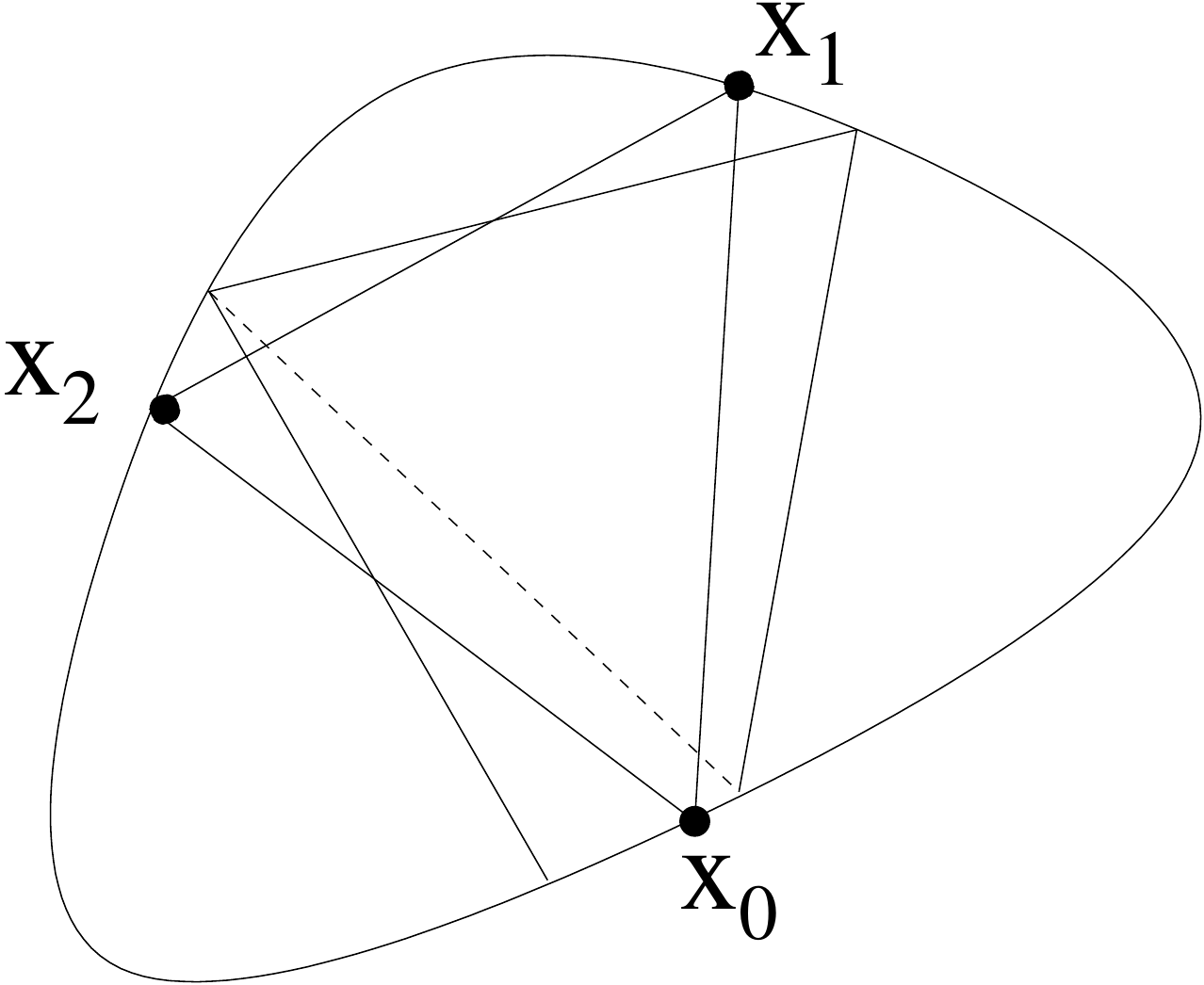}
\vspace{5mm}
\caption{
The periodic orbit $(x_0,x_1,x_2)$ has a well defined 
length function $L$ (the perimeter). A nearby orbit does not have to be
 three-periodic. 
To extend the length function to such orbits, replace the third segment with 
the one connecting the second reflection point with the starting point 
(dashed line). The extended length function is then the perimeter of
 the obtained triangle.}
\label{F:Bill}
\end{center}
\end{figure}

The function $L$ in the above statement is only defined on $P^3$. However, since the billard boundary is smooth, the function $L$ can be extended to a $\mathcal C^1$-differentiable function in a neighborhood of every three-period orbit. Indeed, one can just replace the second reflection with the straight line connecting the second collision point ($x_2$) with the initial point ($x_0$); see Figure~\ref{F:Bill}.

Now, since the billiard boundary is $\mathcal C^3$, both sides of Equation~(\ref{E:Main}) are defined and $\mathcal C^1$-differentiable in an open set containing $P^3$. Note, however, that, in general, Equation~(\ref{E:Main}) does not have to hold in this whole open set.

\subsection{Proof of theorem \ref{T:main}}
Assume to the contrary that ${\rm dim}_{\mathcal H}(P^3) = s$ with $1< s\leq 2$. Note that the set $P^3$ is closed. 
Then, by Theorem \ref{T:SMD}, for any $s^{\prime},\ 1 <s^{\prime} < s$, there exists
 a set $E^{\prime} \subseteq P^3$, such that  $D_{(t,\theta)} T^3 =\id$ for any $(t,\theta)\in E'$. According to 
Proposition \ref{P:2func}, on this set, the identity 
 \[
k(t) L(t,\theta) = 2 \sin^3(\theta)
\]
holds.

We apply Theorem  \ref{T:SDF} to $E^{\prime} \subseteq P^3$ 
and the functions $f(t,\theta) = k(t) L(t,\theta)$, $g(t,\theta)=2 \sin^3(\theta)$.
It gives a full ${\mathcal H}^{s'}$-measure subset 
$E^{\prime \prime}\subseteq E^{\prime}$ such that there are two linearly independent unit 
vectors $\bf{v_1}$ and $\bf{v_2}$ with
$f_{\bf{v_j}}(p)=g_{\bf{v_j}}(p),\ j=1,2$, for any $p\in E^{\prime \prime} $. In particular, since $f$ and $g$ are $\mathcal C^1$-differentiable, 
we have the equality of the partial derivatives: 
\[
\frac{\partial}{\partial \theta } (k(t) L(t,\theta)) = \frac{\partial}{\partial \theta } 
( 2 \sin^3(\theta)).
\]  
This leads to a contradiction since $\partial_{\theta}L(t,\theta) = 0$ by Fermat's principle and 
 $\partial_{\theta}\sin^3(\theta) = 0$ only if $\theta = 0, \pi/2, \pi$. However, there are no three-period 
orbits for these values of $\theta$.
\qed

\subsection{Proof of theorem \ref{T:1reg}}

We may assume that $0<\mathcal H^1(P^3)\leq\infty$. Let $E\subseteq P^3$ be the set that consists of all points $p$ such that $P^3$ does not have a tangent line at $p$. If $\mathcal H^1(E)=0$, we are done. If $\mathcal H^1(E)>0$, the proof follows the lines of the proof of Theorem~\ref{T:main}, where one should replace Theorems~\ref{T:SMD} and~\ref{T:SDF} by Theorems~\ref{T:Diff} and~\ref{T:Dir}, respectively. \qed



\end{document}